\documentclass[a4,12pt]{article}
\usepackage{amsmath,amssymb,amsthm}

\newcommand{\R}{\mathbb{R}}
\newcommand{\Z}{\mathbb{Z}}

\newcommand{\1}{\mathbf{1}}
\newcommand{\lam}{\lambda}

\newcommand{\proj}{\varphi}

\newcommand{\la}{\langle}
\newcommand{\ra}{\rangle}

\newtheorem{thm}{Theorem}[section]
\newtheorem{lmm}[thm]{Lemma}
\newtheorem{prp}[thm]{Proposition}

\theoremstyle{remark}
\newtheorem{rmk}[thm]{Remark}

\title{6-transposition property of $\tau$-involutions of \\ vertex operator algebras}
\author{Shinya Sakuma }
\date{ {\small \it 
 Department of Mathematics, National Cheng Kung University, \\
 Tainan 701, Taiwan} }

\begin{document}
\maketitle

\begin{abstract}
In this paper, we study the subalgebra generated by
 two Ising vectors in the Griess algebra of a vertex operator algebra.
We show that the structure of it 
 is uniquely determined by some inner products of Ising vectors. 
We prove that the order of the product of two $\tau$-involutions is
 less than or equal to $6$ and
 we determine the inner product of two Ising vectors.
\end{abstract}

\section{Introduction}

The Monster simple group $\mathbb{M}$ was first constructed by Griess \cite{G}
 as the automorphism group of a commutative non-associative algebra
 of dimension 196884 with a positive-definite invariant bilinear form,
 which is called the Monstrous Griess algebra.
It is known that $\mathbb{M}$ is a 6-transposition group, 
 that is, $\mathbb{M}$ is generated by some $2A$-involutions and 
 $|\tau\tau'|\leq 6$ for any $2A$-involution $\tau$ and $\tau'$. 
Moreover, the conjugacy class of $\tau\tau'$ is one of the nine classes
 $1A$, $2A$, $3A$, $4A$, $5A$, $6A$, $4B$, $2B$ and $3C$ \cite{C,ATLAS}.
It is shown by Conway \cite{C} that
 each $2A$-involution $\tau$ defines a unique idempotent $e_{\tau}$, 
 called an axis, in the monstrous Griess algebra
 such that the inner product of $e_{\tau}$ and $e_{\tau'}$ 
 is uniquely determined by the conjugacy class of the product $\tau\tau'$.

On the other hand, 
 from the point of view of vertex operator algebras (VOAs),
 the Monster simple group is realized as the automorphism group of 
 the Moonshine VOA $V^{\natural}=\oplus_{n=0}^\infty V^{\natural}_n$,
 which is constructed by Frenkel, Lepowsky, Meurman \cite{FLM} 
 and Miyamoto \cite{M3}.
The weight 2 subspace $V^{\natural}_2$ has a structure
 of a commutative algebra which coincides with the monstrous Griess algebra.
Each axis of the algebra is essentially
 a half of a conformal vector $e$ of $V^{\natural}$ with central charge $1/2$
 which generates the Virasoro VOA $L(1/2,0)$.
Such a vector $e$ is called an Ising vector and 
 defines an involutive automorphism $\tau_e$ of $V^{\natural}$
 called a $\tau$-involution 
 by using the symmetry of the fusion rules for $L(1/2,0)$.
Then, $\tau_e$ is a $2A$-involution of $\mathbb{M}$ and 
 we have a one-to-one correspondence 
 between the 2A-involutions of $\mathbb{M}$ 
 and the Ising vectors of $V^{\natural}$.
It is shown in \cite{C}
 that the structure of the subalgebra generated by two Ising vector $e$ and $f$
 in the algebra $V^{\natural}_2$ depends on only the conjugacy class 
 of $\tau_e\tau_f$ and the inner product $\la e,f \ra$ is given by
 the following table: 
$$
\begin{array}{c|ccccccccc}
(\tau_e\tau_f)^{\mathbb{M}} 
 & 1A & 2A & 3A & 4A & 5A & 6A & 3C & 4B & 2B \\ \hline
 \la e , f \ra & 1/4 & 1/2^5 & 13/2^{10} & 1/2^7 & 3/2^9 & 5/2^{10} &
 1/2^8 & 1/2^8 & 0
\end{array}
$$

In this paper, 
 we consider a VOA $V=\oplus_{n= 0}^\infty V_n$ 
 over the field $\R$ of real numbers
 with a positive-definite invariant bilinear form 
 such that $V_0=\R \1$ and $V_1=0$ as the Moonshine VOA.
Then, defining the product $u\cdot v$ by $u_{(1)}v$ for $u,v\in V_2$,
 the weight 2 subspace $V_2$ has a structure of
 a commutative non-associative algebra
 with the positive-definite symmetric invariant bilinear form
 $\la \cdot , \cdot \ra$ such that $u_{(3)}v=\la u,v \ra \1$.
This algebra is called the (general) Griess algebra of $V$.
An Ising vector $a$ is twice an idempotent of the Griess algebra
 with $\la a,a \ra = 1/4$
 and generates the Virasoro VOA $L(1/2,0)$, which is rational
 and has only three irreducible module 
 $L(1/2,0)$, $L(1/2,1/2)$ and $L(1/2,1/16)$.
For an Ising vector $a\in V$, 
 we have a decomposition $V=V_a(0)\oplus V_a(1/2)\oplus V_a(1/16)$,
 where $V_a(h)$ denote the sum of submodules isomorphic to $L(1/2,h)$
 as a module over the Virasoro VOA generated by $a$.
We define a endomorphism $\tau_a$ of $V$ by
$$
\tau_a = \left\{
\begin{array}{lll}
1 & {\rm on} & V_a(0)\oplus V_a(1/2) \\
-1 & {\rm on} & V_a(1/16).
\end{array} \right.
$$
Then $\tau_a$ is an automorphism of $V$ such that
 $\sigma^{-1}\tau_a \sigma =\tau_{a^\sigma}$ for any $\sigma\in {\rm Aut}(V)$.

The main result of this paper is that 
 $\tau$-involutions of a VOA $V$ satisfy a 6-transposition property,
 that is, $|\tau_e\tau_f|\leq 6$ for any Ising vector $e$ and $f$.
The subgroup $T$ generated by $\tau_e$ and $\tau_f$
 is a Dihedral group such that $(\tau_e\tau_f)^N=1$,
 where $N=|e^T\cup f^T|$ and 
 $a^T$ denotes the orbit of $a\in V$ under the action of $T$. 
Actually, we show that $N\leq 6$.
To prove this, 
 we study the subalgebra $B_{e,f}$ of the Griess algebra $B=V_2$
 generated by two Ising vectors $e$ and $f$.
We show that
 the dimension of $B_{e,f}$ is less than or equal to 8
 and the structure of $B_{e,f}$ is uniquely determined
 by $\la e,f \ra$ and $\la e,e^{\tau_f} \ra$.
Moreover, we determine the inner product of two Ising vectors,
 that is, 
$$
\begin{array}{l}
 \mbox{if $N=2$, then } \la e,f \ra =0 \mbox{ or } 1/2^5, \\
 \mbox{if $N=3$, then } \la e,f \ra =13/2^{10} \mbox{ or } 1/2^8, \\
 \mbox{if $N=4$, then }
 (\la e,f \ra , \la e, e^{\tau_f}\ra )=(1/2^7,0 )
 \mbox{ or } (1/2^8,1/2^5 ), \\
 \mbox{if $N=5$. then } \la e,f \ra=\la e, e^{\tau_f}\ra=3/2^{9} \
 \mbox{ and } \\
 \mbox{if $N=6$, then } \la e,f \ra=5/2^{10}, \
 \la e,e^{\tau_f}\ra=13/2^{10} \mbox{ and } \la e^{\tau_f},f^{\tau_e}\ra=1/2^5.
\end{array}
$$
The inner product $\la e,f \ra$ for the case of $N=2$ and $3$ 
 was shown in \cite{M1,M2}.
We note that 
 these inner products are just given in the case of the Moonshine VOA and
 $B_{e,f}$ coincides with the subalgebra of the Griess algebra
 constructed in \cite{C}, \cite{LYY1} and \cite{LYY2}.

For each Ising vector $a$, we have an eigenspace decomposition of
 the Griess algebra $B$ with respect to the adjoint action of $a$:
 $B=\R a \oplus E_a(0)\oplus E_a(1/2)\oplus E_a(1/16)$,
 where $E_a(h)=\{ v\in B | a\cdot v = h v \}$ for $h=0,1/2,1/16$.
By using this decomposition, we calculate some products in $B$. By
definition of $\tau_a$, we can see that for $x\in B$,
 $x-x^{\tau_a}$ is an eigenvector with eigenvalue $1/16$
 and so $ax-\frac{1}{16}x$ is fixed by $\tau_a$.
For any $x\in B$, the subspace $X_a(x)$
 spanned by $a$, $x$, $x^{\tau_a}$ and $\alpha(a,x):=ax-\frac{1}{16}(a+x)$
 is invariant under the adjoint action of $a$ and
 the projections of $x$ to $E_a(1/2)$ and $E_a(0)$ are expressed
 by $a$, $\frac{1}{2}(x+x^{\tau_a})=:x^+$ and $\alpha(a,x)$ in $X_a(x)$.
Then, by using the eigenspace decomposition of an Ising vector $a$
 on $B$ and fusion rules for $L(1/2,0)$,
 we calculate $\alpha(a,x)\cdot \alpha(a,y)$ in $\sum_u X_a(u)$
 where $u$ runs in $\{ x, y, \alpha(a,x)\cdot y+\alpha(a,y)\cdot x, x^+y^+ \}$.

Let $X(a,b)=X_a(b)+X_a(a^{\tau_b})+X_a(\alpha(b,b^{\tau_a}))$
 for Ising vectors $a$ and $b$.
By using the above method,
 we calculate $\alpha(e,f)\cdot \alpha(e,f)$ in $X(e,f)$ and
 $\alpha(f,e)\cdot \alpha(f,e)$ in $X(f,e)$.
Since $\alpha(e,f)=\alpha(f,e)$ by the definition, we have
 $\alpha(e,f)\cdot \alpha(e,f)=\alpha(f,e)\cdot \alpha(f,e)$
 and this relation implies that $X(e,f)=X(f,e)$.
By using these relations,
 we show that $X(e,f)$ is closed under the multiplication of $B$.
In particular, $X(e,f)=B_{e,f}$.

It is shown in \cite{M3} that
 $0 \leq \la a,b \ra \leq \frac{1}{12}$ for distinct Ising vectors $a$ and $b$.
By using this result and relations in $B_{e,f}$,
 we can show that $|e^T\cup f^T|\leq 6$.
This implies $|\tau_e\tau_f| \leq 6$.

\smallskip

\noindent {\bf Acknowledgement} 
The author wishes to thank 
 Ching Hung Lam and Masahiko Miyamoto
 for stimulating discussions and useful comments. 

\section{Griess Algebra of Vertex Operator Algebras}

Let $V=(V,Y,\1,\omega)$ be a vertex operator algebra (VOA) and
 let $Y(v,z)=\sum_{n\in\Z}v_{(n)}z^{-n-1}$
 denote the vertex operator of $V$ for $v\in V$,
 where $v_{(n)}\in {\rm End}(V)$.
Throughout this paper, we only consider VOAs
 over the field $\R$ of real numbers with the following grading:
$$
V=\oplus_{n=0}^\infty V_n, \ V_0=\R \1 \mbox{ and } V_1=0.
$$
Then, $V$ has a unique invariant bilinear form
 $\la \cdot , \cdot \ra$ with $\la \1 , \1 \ra = 1$.
We also assume that $\la \cdot , \cdot \ra$ is positive-definite.

\subsection{Ising vector and $\tau$-involution}

A vector $e\in V_2$ is called
 a {\it conformal vector with central charge} $c_e$
 if it satisfies $e_{(1)}e=2e$ and $e_{(3)}e=\dfrac{c_e}{2}\1$.
Then the operators $L^e_{(n)}:=e_{(n+1)}, n\in \Z,$ satisfy
 the Virasoro commutation relation
 $$
 [L^e_{(m)},L^e_{(n)}]=(m-n)L^e_{(m+n)}+\delta_{m+n,0}\frac{m^8-m}{12}c_e
 $$
 for $m,n\in \Z$.
A conformal vector $e\in V_2$ with central charge $1/2$
 is called an {\it Ising vector}
 if $e$ generates the Virasoro VOA $L(1/2,0)$,
The Virasoro VOA $L(1/2,0)$ is rational and has only three irreducible modules
 $L(1/2,0)$, $L(1/2,1/2)$ and $L(1/2,1/16)$.
In \cite{DMZ}, the fusion rules of these modules are given by
$$
\begin{array}{rcl}
L\left( \frac{1}{2},\frac{1}{2} \right) \times
 L\left( \frac{1}{2},\frac{1}{2} \right) & = &
 L\left( \frac{1}{2},0 \right) , \\[1mm]
L\left( \frac{1}{2},\frac{1}{2} \right) \times
 L\left( \frac{1}{2},\frac{1}{16} \right) & = &
 L\left( \frac{1}{2},\frac{1}{16} \right), \\[1mm]
L\left( \frac{1}{2},\frac{1}{16} \right) \times
 L\left( \frac{1}{2},\frac{1}{16} \right) & = &
 L\left( \frac{1}{2},0 \right) + L\left( \frac{1}{2},\frac{1}{2} \right) .
\end{array}
$$
For an Ising vector $e\in V$ and $h=0$, $1/2$ and $1/16$,
 let $V_e(h)$ be the sum of submodules isomorphic to $L(1/2,h)$
 as a module over the Virasoro VOA generated by $e$.
Then, we have the decomposition
$$
V = V_e(0)\oplus V_e(1/2)\oplus V_e(1/16).
$$
Define an endomorphism $\tau_e$ on $V$ by
$$
\tau_e = \left\{
\begin{array}{rl}
1 & on \ V_e(0)\oplus V_e(1/2) \\ -1 & on \ V_e(1/16).
\end{array}
\right.
$$
Then, by the fusion rules for $L(1/2,0)$,
 $\tau_e$ is an involutive automorphism of the VOA $V$ such that
 $\sigma^{-1} \tau_e \sigma = \tau_{e^{\sigma}}$
 for any $\sigma\in {\rm Aut}(V)$.
This automorphism is called a {\it $\tau$-involution}.

\subsection{Griess algebra}

We consider the weight two subspace $B:=V_2$.
For $x,y\in B$, we can
define a product $x\cdot y$ by $x_{(1)}y$
 and the bilinear form $\la x , y \ra$ is given by $\la x , y \ra \1 = x_{(3)}y$.
Then,
 $B$ is a commutative non-associative algebra
 with the symmetric bilinear form $\la \cdot ,\cdot \ra$
 which is invariant, that is, satisfies
 $\la x\cdot y, z\ra=\la x,y\cdot z \ra $ for $x,y,z\in B$.
This algebra $B$ is called the {\it Griess algebra} of $V$. Note
that by definition an Ising vector $a$ is twice an idempotent of
 the algebra $B$ with $\la a,a \ra=\frac{1}{4}$.

Let $a$ be an Ising vector of $V$.
Let $B_a(h)=B\cap V_a(h)$ and $E_a(h)=\{ v\in B | a\cdot v=hv \}$
 for an Ising vector $a\in V$ and $h=0$, $1/2$ and $1/16$.
Then, we have an eigenspace decomposition
\begin{equation}\label{eigenspace}
B_a(0)=\R a\oplus E_a(0), B_a(1/2)=E_a(1/2) \mbox{ and }
B_a(1/16)=E_a(1/16)
\end{equation}
(see \cite{M1}). Let $\proj_a^{\pm}(x) = \dfrac{1}{2}(x \pm
x^{\tau_a})$ for $x\in B$. Then, $\proj_a^{+}(x)$ and
$\proj_a^{-}(x)$ are the projections of $x$
 to $B_a^+:=B_a(0) \oplus B_a(1/2)$ and $B_a(1/16)$, respectively.
Thus, by (\ref{eigenspace}),
 $a \cdot \proj_a^{-}(x) = \dfrac{1}{16} \proj_a^{-}(x)$
 and so $ax-\dfrac{1}{16} x\in B_a^+$.
Let $\proj_a^{i}(x)$ be the projection of $x\in B$ to $B_a(i/2)$
 for $i=0,1$ and define $(a|x):=\la a,x \ra / \la a,a \ra = 4 \la a,x \ra $.
Then, we see that
 $ a \cdot x = 2 (a|x) a +\dfrac{1}{2}\proj_a^{1}(x) +\dfrac{1}{16}\proj_a^{-}(x)$
 by (\ref{eigenspace}), and so
\begin{eqnarray}\label{proj1/2}
\proj_a^{1}(x) & = & 2ax-4(a|x)a-\frac{1}{8}\proj_a^{-}(x), \label{proj1}
\end{eqnarray}
and
\begin{eqnarray}
a(ax) & = &
 4 (a|x) a +\frac{1}{4}\proj_a^{1}(x) +\frac{1}{2^8}\proj_a^{-}(x) \nonumber \\
 & = & \frac{1}{2}ax +3(a|x)a -\frac{7}{2^8}\proj_a^{-}(x).
\label{a(ax)}
\end{eqnarray}
Let $\alpha(a,x) = ax-\dfrac{1}{16}(a+x)$ for $x\in B$. Then,
$\alpha(a,x)\in B_a^+$ and by (\ref{a(ax)}) we have
\begin{eqnarray}
a \cdot \alpha(a,x) & = & a(ax)-\frac{1}{16}(aa+ax) \nonumber \\
 & = & \frac{7}{16}ax+\Big( 3(a|x)-\frac{1}{8} \Big) a
  -\frac{7}{2^8}\proj_a^{-}(x) \nonumber \\
 & = & \frac{7}{16}\alpha(a,x) + \Big( 3(a|x)-\frac{25}{2^8} \Big) a
  +\frac{7}{2^8}\proj_a^{+}(x) \label{a alpha(a,x)}
\end{eqnarray}
and
\begin{eqnarray}
\la a , \alpha(a,x) \ra & = & \la a , ax \ra -\frac{1}{16}(\la a , a \ra +\la a , x \ra ) \nonumber \\
 & = & 2\la a , x \ra -\frac{1}{16}\Big( \frac{1}{4} +\la a , x \ra \Big) \nonumber \\
 & = & \frac{31}{16}\la a , x \ra -\frac{1}{2^6}. \label{(a|alpha(a,x))}
\end{eqnarray}

Set
$A_a(x,y)  =   \alpha(a,x) \cdot y + \alpha(a,y) \cdot x$
for an Ising vector $a$ and $x,y\in B$.
\begin{prp}\label{proj_a^1}
Let $a$ be an Ising vector and $x,y\in B$. Then, we have
\begin{eqnarray}
\proj_a^1 (x^+y^+ ) 
 = \frac{8}{3}
 \proj_a^1 \left( A_a(x,y) - \Big( \Big( (a|x)-\frac{1}{2^5} \Big) y
 + \Big((a|y)-\frac{1}{2^5} \Big) x \Big) \right) \label{x^+y^+}
\end{eqnarray}
and
\begin{eqnarray}
\proj_a^1(x)\proj_a^1(y) & = &
 \proj_a^+ ( A_a(x,y)) +\frac{1}{8}x^+y^+
 -\proj_a^1 \left( A_a(x,y) +\frac{1}{8} x^+y^+ \right) \nonumber \\
& & - \left( 8(a|x)(a|y)-\frac{1}{8}((a|x)+(a|y)) \right) a, \label{x^1y^1}
\end{eqnarray}
where $x^+=\proj_a^+(x)$ and $y^+=\proj_a^+(y)$.
\end{prp}

\begin{proof}
To simplify, we sometimes denote $\proj_a^i(z)$ by $z^i$
 for $z\in B$ and $i=0,1,\pm$.
By the fusion rules for $L(1/2,0)$, if $u\in B_a^+$ and $v\in B$
then
\begin{eqnarray*}
\proj_a^1 \left( uv \right) & = & u^0 v^1 + u^1 v^0 \\
\proj_a^0 \left( uv \right) & = & u^0 v^0 + u^1 v^1
\end{eqnarray*}
 and
\begin{eqnarray*}
\proj_a^0 (\alpha(a,v)) & = &
\left( 2(a|v)-\frac{1}{16}\right) a-\frac{1}{16}v^0, \\
\proj_a^1 (\alpha(a,v)) & = & \left( \frac{1}{2}-\frac{1}{16}
\right) v^1 \ = \ \frac{7}{16}v^1
\end{eqnarray*}
since $\alpha(a,v)=av-\frac{1}{16}(a+v)$. Hence,
\begin{eqnarray*}
\proj_a^1(A_a(x,y))
 & = & \proj_a^1(\alpha(a,x) \cdot y )+\proj_a^1(\alpha(a,y) \cdot x ) \\
 & = & \left( \Big( 2(a|x)-\frac{1}{16} \Big) a-\frac{1}{16}x^0 \right) y^1
  +\frac{7}{16}x^1y^0 \\
 &  & +\left( \Big( 2(a|y)-\frac{1}{16} \Big) a-\frac{1}{16}y^0 \right) x^1
  +\frac{7}{16}y^1x^0 \\
 & = & \left( (a|x)-\frac{1}{2^5}\right) y^1
  +\left( (a|y)-\frac{1}{2^5}\right) x^1
  +\frac{3}{8} (x^0y^1 + x^1y^0) .
\end{eqnarray*}
Therefore, this implies (\ref{x^+y^+})
 since $\proj_a^1 (x^+y^+ )=x^0y^1 + x^1y^0$.
Next,
\begin{eqnarray*}
&  & \hspace{-10mm} \proj_a^0 (A_a(x,y)) \\ & = &
 \left( \Big( 2(a|x)-\frac{1}{16} \Big) a-\frac{1}{16}x^0 \right) y^0
  +\frac{7}{16}x^1y^1 \\
 &  & +\left( \Big( 2(a|y)-\frac{1}{16} \Big) a-\frac{1}{16}y^0 \right) x^0
  +\frac{7}{16}y^1x^1 \\
 & = & \left( 2(a|y) \Big( 2(a|x)-\frac{1}{16} \Big)
  +2(a|x) \Big( 2(a|y)-\frac{1}{16} \Big) \right) a
  -\frac{1}{8}x^0y^0+\frac{7}{8}x^1y^1 \\
 & = & \left( 8(a|x)(a|y)-\frac{1}{8}((a|x)+(a|y)) \right) a
  -\frac{1}{8}\proj_a^0(x^+y^+) + x^1y^1,
\end{eqnarray*}
by $x^0y^0 = \proj_a^0(x^+y^+) - x^1y^1$.
Thus,
\begin{eqnarray*}
& & \hspace{-10mm} x^1y^1 + \left( 8(a|x)(a|y)-\frac{1}{8}((a|x)+(a|y)) \right) a \\
& = & \proj_a^0 \left( (A_a(x,y) +\frac{1}{8}x^+y^+ \right) \\
& = & \proj_a^+ ( A_a(x,y)) +\frac{1}{8} x^+y^+
 -\proj_a^1 \left( A_a(x,y) +\frac{1}{8} x^+y^+ \right) ,
\end{eqnarray*}
which implies (\ref{x^1y^1}).
\end{proof}

\begin{prp}\label{prp alpha(a,x)alpha(a,y)}
For an Ising vector $a$ and $x,y\in B$, we have
\begin{eqnarray}
&  & \hspace{-10mm} \alpha(a,x) \cdot \alpha(a,y) \nonumber \\
& = & \frac{3}{16}\proj_a^+ (A_a(x,y)) +\frac{1}{2^8}x^+y^+
 +\Big( 6(a|x)(a|y)-\frac{7}{2^5}((a|x)+(a|y)) +\frac{1}{2^7} \Big) a
 \nonumber \\
 &  & -\frac{1}{3}\proj_a^1 \left( A_a(x,y)
  -\frac{7}{4} \Big( \Big( (a|x)-\frac{1}{2^5} \Big) y
   + \Big( (a|y)-\frac{1}{2^5} \Big) x \Big) \right) .
\end{eqnarray}
\end{prp}

\begin{proof}
For $v\in B$, since
$\alpha(a,v)=\proj_a^+(\alpha(a,v))=\alpha(a,v^+)$,
$$\alpha(a,v) +\frac{1}{16}v^+ = a \cdot v^+ -\frac{1}{16}a
= \Big( 2(a|v)-\frac{1}{16} \Big) a +\frac{1}{2}v^1. $$ Thus, by
using Proposition \ref{proj_a^1},
\begin{eqnarray*}
&  & \hspace{-.3in} \left( \alpha(a,x) +\frac{1}{16}x^+ \right)
 \left( \alpha(a,y) +\frac{1}{16}y^+ \right) \\
 & = & \left( \Big( 2(a|x)-\frac{1}{16} \Big) a +\frac{1}{2}x^1 \right)
 \left( \Big( 2(a|y)-\frac{1}{16} \Big) a +\frac{1}{2}y^1 \right) \\
 & = & \frac{1}{4} x^1y^1
 +2 \Big( 2(a|x)-\frac{1}{16} \Big)\Big( 2(a|y)-\frac{1}{16} \Big) a \\
 &  & +\frac{1}{4} \left( \Big( 2(a|x)-\frac{1}{16} \Big) y^1
 +\Big( 2(a|y)-\frac{1}{16} \Big) x^1 \right) \\
 & = & \frac{1}{4} \proj_a^+ ( A_a(x,y)) +\frac{1}{2^5} x^+y^+
 -\frac{1}{4} \proj_a^1 \left( A_a(x,y) +\frac{1}{8} x^+y^+ \right) \\
 &  & - \left( 2(a|x)(a|y)-\frac{1}{2^5}((a|x)+(a|y))
 -2 \Big( 2(a|x)-\frac{1}{16} \Big)\Big( 2(a|y)-\frac{1}{16} \Big) \right) a \\
 &  & +\frac{1}{2} \left( \Big( (a|x)-\frac{1}{2^5} \Big) y^1
 +\Big( (a|y)-\frac{1}{2^5} \Big) x^1 \right) \\
 & = &  \frac{1}{4} \proj_a^+ ( A_a(x,y)) +\frac{1}{2^5} x^+y^+
 +(6(a|x)(a|y)-\frac{7}{2^5}((a|x)+(a|y)) +\frac{1}{2^7})a \\
 &  & -\frac{1}{3}\proj_a^1 \left( A_a(x,y)
  -\frac{7}{4} \Big( \Big( (a|x)-\frac{1}{2^5} \Big) y
   + \Big( (a|y)-\frac{1}{2^5} \Big) x \Big) \right) .
\end{eqnarray*}
On the other hand,
\begin{eqnarray*}
&  & \hspace{-.3in} \left( \alpha(a,x) +\frac{1}{16}x^+ \right)
 \left( \alpha(a,y) +\frac{1}{16}y^+ \right) \\
 & = & \alpha(a,x) \alpha(a,y)
  +\frac{1}{16}(\alpha(a,x) y^+ +\alpha(a,y) x^+ ) +\frac{1}{2^8}x^+y^+ \\
 & = & \alpha(a,x) \alpha(a,y) +\frac{1}{16}\proj_a^+ (A_a(x,y))
  +\frac{1}{2^8}x^+y^+ .
\end{eqnarray*}
Therefore, we get
\begin{eqnarray*}
&  & \hspace{-.3in} \alpha(a,x) \alpha(a,y) \\
& = & \frac{3}{16}\proj_a^+ (A_a(x,y)) +\frac{1}{2^8}x^+y^+
 +(6(a|x)(a|y)-\frac{7}{2^5}((a|x)+(a|y)) +\frac{1}{2^7})a \\
 &  & -\frac{1}{3}\proj_a^1 \left( A_a(x,y)
  -\frac{7}{4} \Big( \Big( (a|x)-\frac{1}{2^5} \Big) y
   + \Big( (a|y)-\frac{1}{2^5} \Big) x \Big) \right) .
\end{eqnarray*}
\end{proof}

\section{Subalgebra Generated by Two Ising Vectors}

In this section, we study the subalgebra $B_{e,f}$ of $B$
 generated by two Ising vector $e$ and $f$.
We show that the structure of $B_{e,f}$ is uniquely determined
 by $(e|f)$ and $(e|e^{\tau_f})$.

For any Ising vectors $a,b\in V$, let $T_{a,b}=\la \tau_a , \tau_b
\ra$ denote the subgroup of $Aut(V)$ generated by $\tau_a$ and
$\tau_b$.
By definition of $\alpha(a,b)$,
 we see that $\alpha(a,b)=\alpha(b,a)$ and $\alpha(a,b)$ is fixed by $T_{a,b}$.
\begin{prp}\label{inner product}
Let $e$ and $f$ be Ising vectors and $\rho=\tau_e\tau_f$. Then, we
have $\la e , e^{\rho^n} \ra = \la f , f^{\rho^n} \ra$ for any $n\in
\Z$.
\end{prp}

\begin{proof}
Set $\alpha=\alpha(e,f)$. Then by (\ref{a alpha(a,x)}),
\begin{eqnarray*}
e \cdot \alpha & = & \frac{7}{16}\alpha + \Big( 3(e|f)-\frac{25}{2^8} \Big) e
 +\frac{7}{2^9}(f^{\tau_e}+f)
\end{eqnarray*}
and so
\begin{eqnarray*}
\la f^{\rho^n}e , \alpha \ra
 & = & \frac{7}{16}\la f^{\rho^n} , \alpha \ra
 + \Big( 3(e|f)-\frac{25}{2^8} \Big) \la f^{\rho^n} , e \ra \\
& & \hspace{5mm}
 +\frac{7}{2^9}(\la f^{\rho^n} , f^{\tau_e} \ra +\la f^{\rho^n} , f \ra ).
\end{eqnarray*}
By exchanging $e$ and $f$ in the above equation,  since
$\alpha(e,f)=\alpha(f,e)$, we obtain
\begin{eqnarray*}
\la e^{\rho^{-n}}f , \alpha \ra
 & = & \frac{7}{16}\la e^{\rho^{-n}} , \alpha \ra
 + \Big( 3(f|e)-\frac{25}{2^8} \Big) \la e^{\rho^{-n}} , f \ra \\
& & \hspace{5mm}
 +\frac{7}{2^9}(\la e^{\rho^{-n}} , e^{\tau_f} \ra +\la e^{\rho^{-n}} , e \ra ).
\end{eqnarray*}
On the other hand, since $\alpha$ is fixed by $T_{e,f}$, we see that
\begin{eqnarray*}
\la f^{\rho^n}e , \alpha \ra
 & = & \la fe^{\rho^{-n}} , \alpha^{\rho^{-n}} \ra
 \ = \ \la e^{\rho^{-n}}f , \alpha \ra, \\
\la f^{\rho^n} , e \ra & = & \la f , e^{\rho^{-n}} \ra \ = \ \la
e^{\rho^{-n}} , f \ra
\end{eqnarray*}
and by (\ref{(a|alpha(a,x))})
\begin{eqnarray*}
\la f^{\rho^n} , \alpha \ra & = & \la f , \alpha \ra
 \ = \ \frac{31}{16}\la e , f \ra -\frac{1}{2^6} \\
 & = & \la e , \alpha \ra
 \ = \ \la e^{\rho^{-n}} , \alpha \ra
\end{eqnarray*}
for $n\in\Z$. Hence, by using these equations, we have
$$
\la f^{\rho^n} , f^{\tau_e} \ra +\la f^{\rho^n} , f \ra = \la
e^{\rho^{-n}} , e^{\tau_f} \ra +\la e^{\rho^{-n}} , e \ra.
$$
Since $f^{\tau_e}=f^{\rho^{-1}}$, $e^{\tau_f}=e^{\rho}$
 and $\la e^{\rho^{-n}} , e \ra =\la e , e^{\rho^{n}} \ra$,
$$
\la f , f^{\rho^{n+1}} \ra +\la f , f^{\rho^n} \ra = \la e ,
e^{\rho^{n+1}} \ra +\la e , e^{\rho^{n}} \ra
$$
for $n\in\Z$. If $n=0$, then $\la e^{\rho^{n}} , e
\ra=\frac{1}{4}=\la f^{\rho^n} , f \ra$. Therefore, by induction on
$n$, we have $ \la e , e^{\rho^{n}} \ra = \la f , f^{\rho^n} \ra $
for any $n\in\Z$.
\end{proof}

\begin{prp}\label{prp alpha(e,f)alpha(e,f)}
We have
\begin{eqnarray}
&  & \hspace{-.3in} \alpha(e,f) \cdot \alpha(e,f) \nonumber \\
& = & \frac{7}{3} \Big( 4\lam_1^2 -\frac{1}{2^4}\lam_1 -\frac{1}{2^{12}}
 +\frac{1}{2^{6}}\lam_2 \Big) e
 +\frac{7^2}{3\cdot 2^4} \Big( \lam_1 -\frac{5}{2^{8}} \Big) \proj_e^+(f)
 +\frac{7^2}{3\cdot 2^{12}}\proj_e^+(e^{\tau_f}) \nonumber \\
 &  & -\frac{1}{3} \Big( 5\lam_1 + \frac{13}{2^7} \Big) \alpha(e,f)
 -\frac{7}{3\cdot 2^7}\alpha(e,e^{\tau_f}) +\frac{7}{2^9}
 \alpha(f,f^{\tau_e}) , \label{alpha(e,f)alpha(e,f)}
\end{eqnarray}
where $\lam_1=(e|f)$ and $\lam_2=(e|e^{\tau_f})$
\end{prp}

\begin{proof}
Let $x^i=\proj_e^i(x)$ for $x\in B$ and $i=0,1,\pm$.
By Proposition \ref{prp alpha(a,x)alpha(a,y)},
\begin{eqnarray*}
&  & \hspace{-.3in} \alpha(e,f) \cdot \alpha(e,f) \\
& = & \frac{3}{16}\proj_a^+ (A_e(f,f)) +\frac{7}{2^8}f^+f^+
 + \Big( 6\lam_1^2-\frac{7}{2^4}\lam_1 +\frac{1}{2^7} \Big) e \\
 &  & -\frac{1}{3}\proj_a^1 \left( A_e(f,f)
 -\frac{7}{2} \Big( \lam_1-\frac{1}{2^5} \Big) f \right) .
\end{eqnarray*}
Since $f^+=\frac{1}{2}(f+f^{\tau_e})$ by definition,
\begin{equation}\label{f^+f^+} f^+f^+=\frac{1}{2}(f + f^{\tau_e}
+f\cdot f^{\tau_e}) =\frac{17}{16}f^+
+\frac{1}{2}\alpha(f,f^{\tau_e}).
\end{equation}
On the other hand,
 by (\ref{a alpha(a,x)}) and $\alpha(e,f)=\alpha(f,e)$,
\begin{eqnarray}
A_e(f,f) & = & 2f \cdot \alpha(f,e) \nonumber \\
 & = & \frac{7}{8}\alpha(f,e) + \Big( 6\lam_1-\frac{25}{2^7} \Big) f
 +\frac{7}{2^8}(e+e^{\tau_f}) \label{A(f,f)}
\end{eqnarray}
and so, by $\proj_a^1 (\alpha(e,f))=\dfrac{7}{16}f^1$,
\begin{eqnarray*}
& & \hspace{-.3in} \proj_a^1 \left( A_e(f,f)
 -\frac{7}{2} \Big( \lam_1-\frac{1}{2^5} \Big) f \right) \\
 & = & \frac{7^2}{2^7}f^1 + \Big( 6\lam_1-\frac{25}{2^7} \Big) f^1
 +\frac{7}{2^8}(e^{\tau_f})^1
 -\frac{7}{2} \Big( \lam_1-\frac{1}{2^5} \Big) f^1\\
 & = &  \Big( \frac{5}{2}\lam_1+\frac{19}{2^6} \Big) f^1
  +\frac{7}{2^8}(e^{\tau_f})^1.
\end{eqnarray*}
By (\ref{proj1}), for $x\in B$,
\begin{equation}\label{proj_e^1 alpha}
x^1=2ex-4(e|x)e-\frac{1}{8}x^{-} =2\alpha(e,x)- \Big(
4(e|x)-\frac{1}{8} \Big) e +\frac{1}{8}x^{+}.
\end{equation}
Therefore,
\begin{eqnarray*}
&  & \hspace{-.3in} \alpha(e,f) \cdot \alpha(e,f) \\
& = & \frac{3}{16} \left( \frac{7}{8}\alpha(e,f)
 + \Big( 6\lam_1 -\frac{25}{2^7} \Big) f^+
 +\frac{7}{2^8}(e+(e^{\tau_f})^+) \right) \\
 &  & +\frac{7}{2^8} \left( \frac{17}{16}f^+
 +\frac{1}{2}\alpha(f,f^{\tau_e}) \right)
 + \Big( 6\lam_1^2-\frac{7}{2^4}\lam_1 +\frac{1}{2^7} \Big) e \\
 &  & -\frac{1}{3} \Big( \frac{5}{2}\lam_1+\frac{19}{2^6} \Big)
  \left( 2\alpha(e,f) +\frac{1}{8}f^+
   - \Big( 4\lam_1-\frac{1}{8} \Big) e \right) \\
 &  & -\frac{7}{3 \cdot 2^8}
  \left( 2\alpha(e,e^{\tau_f})+\frac{1}{8}(e^{\tau_f})^+
   - \Big( 4\lam_2-\frac{1}{8} \Big) e \right) \\
 & = & \frac{7}{3} \Big( 4\lam_1^2 -\frac{1}{2^4}\lam_1 -\frac{1}{2^{12}}
 +\frac{1}{2^{6}}\lam_2  \Big) e
 +\frac{7^2}{3 \cdot 2^4} \Big( \lam_1-\frac{5}{2^8} \Big) f^+
 +\frac{7^2}{3 \cdot 2^{12}}(e^{\tau_f})^+ \\
 &  & -\frac{1}{3} \Big( 5\lam_1+\frac{13}{2^7} \Big) \alpha(e,f)
 -\frac{7}{3 \cdot 2^8}\alpha(e,e^{\tau_f}) +\frac{7}{2^9}\alpha(f,f^{\tau_e}).
\end{eqnarray*}
\end{proof}

\begin{prp}\label{e-f}
We have
\begin{eqnarray}
& & \hspace{-15mm} \frac{1}{7} \Big( 2^{11}\lam_1^2 -9 \cdot
2^4\lam_1 +\frac{33}{2^{4}} +2^3\lam_2 \Big) (e-f)
 +\Big( 2^4 \lam_1-\frac{3}{8} \Big) (f^{\tau_e}-e^{\tau_f}) \nonumber \\
 & & +\frac{1}{2^4}(e^{\tau_f\tau_e}-f^{\tau_e\tau_f})
 -(\alpha(e,e^{\tau_f}) -\alpha(f,f^{\tau_e}))
\ = \ 0 \label{e-f}
\end{eqnarray}
\end{prp}

\begin{proof}
We can get $\alpha(f,e) \cdot \alpha(f,e)$
 by exchanging $e$ and $f$ in (\ref{alpha(e,f)alpha(e,f)}).
Since $\alpha(e,f)=\alpha(f,e)$,
\begin{eqnarray*}
0 & = & \alpha(e,f) \cdot \alpha(e,f) - \alpha(f,e) \cdot \alpha(f,e) \\
 & = & \frac{7}{3} \Big( 4\lam_1^2 -\frac{1}{2^4}\lam_1 -\frac{1}{2^{12}}
 +\frac{1}{2^{6}}\lam_2 \Big) (e-f) \\
 & & +\frac{7^2}{3\cdot 2^4} \Big( \lam_1 -\frac{5}{2^{8}} \Big)
  (\proj_e^+(f) -\proj_f^+(e) ) +\frac{7^2}{3\cdot 2^{12}}
  ( \proj_e^+(e^{\tau_f}) -\proj_f^+(f^{\tau_e}) ) \\
 &  & - \Big( \frac{7}{3\cdot 2^7}+\frac{7}{2^9} \Big)
 (\alpha(e,e^{\tau_f}) - \alpha(f,f^{\tau_e}) ) \\
 & = & \frac{7}{3} \Big( 4\lam_1^2 -\frac{9}{2^5}\lam_1 -\frac{33}{2^{13}}
 +\frac{1}{2^{6}}\lam_2 \Big) (e-f) \\
 & & +\frac{7^2}{3\cdot 2^4} \Big( \lam_1 -\frac{3}{2^7} \Big)
 (f^{\tau_e}-e^{\tau_f} )
 +\frac{7^2}{3\cdot 2^{12}} ( e^{\tau_f\tau_e}-f^{\tau_e\tau_f} ) \\
 & & -\frac{7^2}{3\cdot 2^9} (\alpha(e,e^{\tau_f}) - \alpha(f,f^{\tau_e}) ).
\end{eqnarray*}
Dividing this by $\frac{7^2}{3\cdot 2^9}$, we obtain the desired equation.
\end{proof}

\begin{prp}\label{prp proj_f^-}
We have
\begin{eqnarray}
& & \hspace{-15mm} \frac{1}{7}\left( 2^{15}\lam_1^2 -2^{9}\lam_1
+2^{7}\lam_2 -9 \right)
 \proj_f^-(e)
 \nonumber \\ & &
 + \Big( 2^{8}\lam_1 -5 \Big) \proj_f^-(f^{\tau_e})
 +\proj_f^-(e^{\tau_f\tau_e})
 \ = \ 0. \label{proj_f^-}
\end{eqnarray}
In particular, $\{ \proj_f^{-}(e), \proj_f^{-}(f^{\tau_e}),
\proj_f^{-}(e^{\tau_f\tau_e}) \}$ is linearly dependent.
\end{prp}

\begin{proof}
Since $f$, $\alpha(e,f)$, $\alpha(e,e^{\tau_f})$ and
$\alpha(f,f^{\tau_e})$
 are fixed by $\tau_f$, and $\proj_f^-(x^{\tau_f})=-\proj_f^-(x)$,
 applying $\proj_f^-(\cdot )$ to (\ref{alpha(e,f)alpha(e,f)}),
\begin{eqnarray*}
0 & = & \frac{7}{3} \Big( 4\lam_1^2 -\frac{1}{2^4}\lam_1
-\frac{1}{2^{12}}
 +\frac{1}{2^{6}}\lam_2 \Big) \proj_f^-(e) \\
 & & +\frac{7^2}{3\cdot 2^5} \Big( \lam_1 -\frac{5}{2^{8}} \Big)
  \proj_f^-(f^{\tau_e})
 +\frac{7^2}{3\cdot 2^{13}}\proj_f^-(e^{\tau_f}+e^{\tau_f\tau_e}) \\
  & = & \frac{7}{3} \Big( 4\lam_1^2 -\frac{1}{2^4}\lam_1
 -\frac{9}{2^{13}}+\frac{1}{2^6}\lam_2 \Big) \proj_f^-(e) \\
 & & +\frac{7^2}{3\cdot 2^5} \Big( \lam_1 -\frac{5}{2^{8}} \Big)
 \proj_f^-(f^{\tau_e}) +\frac{7^2}{3\cdot 2^{13}}\proj_f^-(e^{\tau_f\tau_e}
\end{eqnarray*}
Dividing this by $\frac{7^2}{3\cdot 2^{13}}$, we obtain the desired
equation.
\end{proof}

For an Ising vector $a$ and $x\in B$,
 let $X_a(x)$ be the subspace of $B$ spanned by $a, x, x^{\tau_a}$ and $\alpha(a,x)$.
Let $e$ and $f$ be Ising vectors in $B$ and set
 $X=X_e(f)+X_e(e^{\tau_f})+\R \alpha(f,f^{\tau_e})$.
We consider the following condition $(*)$ for $(u,v)\in X\times X$:
$$
(*) \ \begin{array}{l}
 u\cdot v \in X, \mbox{ and } \\
 u\cdot v \mbox{ and } \la u , v \ra \mbox{ are uniquely determined by }
 \lam_1 \mbox{ and } \lam_2,
\end{array}
$$
where $\lam_1=(e|f)$ and $\lam_2=(e|e^{\tau_f})=(f|f^{\tau_e})$.
Denote by $x^T$ the orbit of $x\in B$ under the action of
 $T=T_{e,f}$. Then,
\begin{lmm}\label{(*)}
We have \\
{\rm (1)} $(a,v)$ satisfies $(*)$ for any $a\in e^T\cup f^T$ and $v\in X$. \\
{\rm (2)} $(\alpha(e,a),\alpha(e,b))$ satisfies $(*)$ for $a,b\in e^T\cup f^T$. \\
{\rm (3)} $(u,v)$ satisfies $(*)$ for any $u,v\in X$. In particular,
$X$ is subalgebra of $B$.
\end{lmm}

\begin{proof}
(1) Since $x^{\tau_e}=x-2\proj_e^-(x)$,
 $e \cdot x^{\tau_e}=e \cdot x -\frac{1}{16}(x-x^{\tau_e})$ for $x\in B$.
By (\ref{a alpha(a,x)}) and (\ref{(a|alpha(a,x))}),
 $e\cdot \alpha(e,x)$ and $\la e,\alpha(e,x) \ra$ are determined by $(e|x)$ for $x\in B$.
Thus, $(e,v)$ satisfies $(*)$ if $v\in X_e(x)$ and $x\in \{ f,e^{\tau_f} \} $.
Let $f^i=\proj_e^i(f)$ for $i=0,1,\pm$. Then, by (\ref{f^+f^+}),
\begin{eqnarray*}
e \cdot \alpha(f,f^{\tau_e}) & = &
 2 (e|\alpha(f,f^{\tau_e}))e +\frac{1}{2} \proj_e^1(\alpha(f,f^{\tau_e})) \\
 & = & 8\la e , \alpha(f,f^{\tau_e}) \ra e
  +\proj_e^1 \left( f^+f^+ -\frac{17}{16}f \right).
\end{eqnarray*}
By Proposition \ref{proj_a^1} and (\ref{A(f,f)}),
\begin{eqnarray*}
\proj_e^1 \left( f^+f^+ -\frac{17}{16}f \right) & = &
 \frac{8}{3}\proj_e^1 \left( \Big( 4\lam_1 -\frac{19}{2^7} \Big) f +\frac{7}{2^8} e^{\tau_f} \right)
\end{eqnarray*}
and by $e \cdot f^{\tau_e}=e \cdot f-\frac{1}{16}(f-f^{\tau_e})$ and
 $\la f , e \cdot f\ra =\la e , f \cdot f\ra =2\la e , f \ra $,
\begin{eqnarray}
\la e , \alpha(f,f^{\tau_e}) \ra & = & \la e , f \cdot f^{\tau_e} \ra -\frac{1}{8} \la e , f^+ \ra \\
 & = & \la f , e \cdot f^{\tau_e} \ra -\frac{1}{8} \la e , f \ra \nonumber \\
 & = & \frac{15}{8}\la e , f \ra -\frac{1}{16} \Big( \frac{1}{4}-\la f , f^{\tau_e} \ra
 \Big) \label{(e|alpha(f,f^{tau_e}))}.
\end{eqnarray}
Since $\proj_e^1 (x)=2\alpha(e,x)- \Big( 4(e|x)-\dfrac{1}{8} \Big) e
 +\dfrac{1}{8}e^{+}(x)$ by (\ref{proj_e^1 alpha}),
 $e \cdot \alpha(f,f^{\tau_e}) $ and $\la e , \alpha(f,f^{\tau_e}) \ra $
 are uniquely determined by $\lam_1$ and $\lam_2$ in $X$.
Hence, $(e,v)$ satisfies $(*)$ for any $v\in X$.
We see that
 $X=X_f(e)+X_f(f^{\tau_e})+\R \alpha(e,e^{\tau_f})$
 by Proposition \ref{e-f}.
Thus, by a similar argument, $(f,v)$ also satisfies $(*)$
 for any $v\in X$.
Since $u^\sigma \cdot v^\sigma =(u \cdot v )^\sigma $ and
 $\la u^\sigma , v^\sigma \ra = \la u,v \ra$
 for $u,v\in X$ and $\sigma\in T$,
 $(a,v)$ satisfies $(*)$ for any $a\in e^T\cup f^T$ and $v\in X$.  \\
(2)
By (1), $\alpha(e,a),\alpha(e,b)\in X$ and
 $A_e(a,b)=\alpha(e,a)\cdot b +\alpha(e,b)\cdot a$
 is determined by $\lam_1$ and $\lam_2$ in $X$.
Thus, so is $\alpha(e,a)\cdot \alpha(e,b)$
 by Proposition \ref{prp alpha(a,x)alpha(a,y)}.
Since
\begin{eqnarray*}
\la \alpha(e,a) , \alpha(e,b) \ra & = & \la a , e\cdot \alpha(e,b)
\ra -\frac{1}{16} \la e+a , \alpha(e,b) \ra,
\end{eqnarray*}
$\la \alpha(e,a) , \alpha(e,b) \ra $ is determined by $\lam_1$ and
$\lam_2$ by (1). Therefore, $(\alpha(e,a),\alpha(e,b))$ satisfies $(*)$. \\
(3) By Proposition \ref{e-f}, we see that $X$ is spanned by
$e,e^{\tau_f},e^{\tau_f\tau_e},f,f^{\tau_e},f^{\tau_e\tau_f},\alpha(e,f)$
and $\alpha(e,e^{\tau_f})$. Therefore, by (1) and (2), $(u,v)$
satisfies $(*)$ for any $u,v\in X$.
\end{proof}
By the above Lemma, we have the following Theorem.
\begin{thm}
Let $B_{e,f}$ denote the subalgebra of $B$ generated by
 two Ising vectors $e$ and $f$.
Then, $\dim B_{e,f} \leq 8$ and
 the structure of $B_{e,f}$ is uniquely determined by
 $(e|f)$ and $(e|e^{\tau_f})$.
\end{thm}

\section{6-transposition Property}

In this section, we prove a 6-transposition property of
 $\tau$-involutions of a VOA.

\begin{lmm}\label{orbit}
Let $e$ and $f$ be Ising vectors and set $T=T_{e,f}$ and $\rho=\tau_e\tau_f$. Then, \\
{\rm (1)} $|e^T|=|f^T|$. In particular, $e=e^{\rho^n}$ if and only if $f=f^{\rho^n}$. \\
{\rm (2)} $e^T= f^T$ if and only if $|e^T|$ is odd and
$f=e^{\rho^{\frac{n+1}{2}}}$, where $n=|e^T|$. \\
{\rm (3)} $(\tau_e\tau_f)^{|e^T\cup f^T|}=1$ as an automorphism of
$V$.
\end{lmm}

\begin{proof}
(1)
Let $a_{2i}=f^{\rho^{i}}$ and $a_{2i-1}=e^{\rho^{i}}$ for $i\in\Z$.
It is easy to see that
 $e^T=\{ a_{2i+1} | i\in\Z\}$ and $f^T=\{ a_{2i} | i\in \Z\}$
 since $T=\la \tau_e, \rho \ra=\la \tau_f, \rho \ra$.
By Proposition \ref{inner product} and $(a_j)^{\rho^i}=a_{2i+j}$
 for $i,j\in\Z$,
 we have $\la a_j,a_k\ra=\la f,a_{k-j}\ra=\la e,a_{k-j-1}\ra$ for any $j,k\in\Z$.
Since $a_j=a_k$ if and only if $\la a_j,a_k\ra=\frac{1}{4}$,
 letting $n=|e^T|$,
 $e^T=\{ a_1, a_3, \dots , a_{2n-1} \}$ and $f^T=\{ a_0, a_2, \dots , a_{2n-2}
 \}$.
Hence, $|e^T|=|f^T|$. \\
(2) Since $\rho^i \tau_f=\tau_f \rho^{-i}$ for $i\in\Z$,
 we see that $a_j^{\tau_f}=a_{-j}=a_{2n-j}$ for $j\in\Z$.
Thus, if $a_j$ in $e^T=\{ a_1, a_3, \dots , a_{2n-1} \}$
 is fixed by $\tau_f$, then $j=n$.
Hence, if $f\in e^T$, then $|e^T|=n$ is odd and
 $f=a_n=e^{\rho^{\frac{n+1}{2}}}$. \\
(3) Let $N=|e^T\cup f^T|$. If $N$ is odd,
 then $e^T=f^T$ and $a_N=e^{\rho^{\frac{N+1}{2}}}=f$. Thus,
$$
\tau_f=\rho^{-\frac{N+1}{2}}\tau_e\rho^{\frac{N+1}{2}}
 =\tau_e\rho^{N+1}=\tau_f\rho^{N}.
$$
If $N$ is even, then $e^T\not=f^T$ and $|e^T|=|f^T|=\frac{N}{2}$.
Then, $a_N=f^{\rho^{\frac{N}{2}}}=f$. Thus,
$$
\tau_f=\rho^{-\frac{N}{2}}\tau_f\rho^{\frac{N}{2}}=\tau_f\rho^{N}.
$$
Therefore, $\rho^{N}=1$.
\end{proof}

\begin{lmm}\label{independent}
If $|e^T\cup f^T|\geq 7$, then $\proj_f^{-}(e)$,
 $\proj_f^{-}(f^{\tau_e})$ and $\proj_f^{-}(e^{\tau_f\tau_e})$
 are linearly independent.
\end{lmm}

\begin{proof}
As in the proof of the above lemma,
 let $a_{2i}=f^{\rho^{i}}$ and $a_{2i-1}=e^{\rho^{i}}$ for $i\in\Z$
 and let $N=|e^T\cup f^T|$.
Set $\lam_j=(f|a_j)=4\la f,a_j \ra$ and
 $\mu_{j,k}=2(\proj_f^{-}(a_j)|\proj_f^{-}(a_k))$ for $j,k\in\Z$.
Consider the $3\times 3$ symmetric matrix $A=(\mu_{j,k})_{1\leq
j,k\leq 3}$. We will show that $\det A > 0$.

It is easy to see that
 $\lam_{k-j}=(a_j|a_k)=(a_{-j}|a_{-k})$
 by $\la a_j, a_k \ra=\la f , a_{j-k} \ra$ and $a_j^{\tau_f}=a_{-j}$
 for $j,k\in\Z$.
Since $\proj_f^{-}(a_j)=\frac{1}{2}(a_j-a_{-j})$,
 $\mu_{j,k}=\lam_{k-j}-\lam_{k+j}$ for $j,k\in \Z$.
As we see in the proof of the above lemma,
 if $0<m<N$, then $f\not=a_m$ and so
 $0\leq \lam_m=4\la f , a_m \ra \leq \frac{1}{3}$
 by Theorem 9.1 of \cite{M3}.
Since $N\geq 7$ by the assumption,
 \begin{equation}\label{mu}
 \begin{array}{ccl}
 \mu_{i,i} = 1-\lam_{2i}\geq \dfrac{2}{3} & \mbox{ for } & i=1,2,3, \\
 \dfrac{1}{3} \geq |\lam_{k-j}-\lam_{k+j}|=|\mu_{j,k}|
 & \mbox{ if } & 1\leq j < k \leq 3 .
 \end{array}
 \end{equation}
Set
 $M_1 = \dfrac{2^3}{3^3}-\dfrac{2}{3}(\mu_{1,2}^2+\mu_{2,3}^2+\mu_{1,3}^2)
 +2\mu_{1,2}\mu_{2,3}\mu_{1,3}$ and
 $M_2 = \det A - M_1 $.
Then, by (\ref{mu}),
 $$M_1 \geq \frac{2^3}{3^3} -\frac{2}{3}\cdot \frac{3}{3^2}
 -2\cdot \frac{1}{3^3}=0$$
and
\begin{eqnarray*}
M_2 & = & \det A - M_1 \\
 & = & \mu_{1,1}\mu_{2,2}\mu_{3,3}
 -(\mu_{1,1}\mu_{2,3}^2+\mu_{2,2}\mu_{1,3}^2+\mu_{3,3}\mu_{1,2}^2) \\
 & & +\frac{2}{3}(\mu_{1,2}^2+\mu_{2,3}^2+\mu_{1,3}^2)-\frac{2^3}{3^3} \\
 & = & \Big( \mu_{1,1}-\frac{2}{3} \Big)(\mu_{2,2}\mu_{3,3}-\mu_{2,3}^2)
 + \Big( \mu_{2,2}-\frac{2}{3} \Big) \Big( \frac{2}{3}\mu_{3,3}-\mu_{1,3}^2 \Big)
 \\ &  &
 + \Big( \mu_{3,3}-\frac{2}{3} \Big) \Big( \frac{2^2}{3^2}-\mu_{1,2}^2 \Big).
\end{eqnarray*}
Since $\mu_{2,2}\mu_{3,3}-\mu_{2,3}^2$,
 $\dfrac{2}{3}\mu_{3,3}-\mu_{1,3}^2$ and
 $\dfrac{2^2}{3^2}-\mu_{1,2}^2$ are positive by (\ref{mu}),
 $M_2\geq 0$.
Thus, $\det A =M_1+M_2 \geq 0$. Assume that $M_2 = 0$. Then,
$\mu_{i,i}=\frac{2}{3}$ for $i=1,2,3$,
 that is, $\lam_2=\lam_4=\lam_6=\frac{1}{3}$.
Thus, $\mu_{1,3}=\lam_2-\lam_4=0$ and so
 $M_1=\dfrac{2^3}{3^3}-\dfrac{2}{3}(\mu_{1,2}^2+\mu_{2,3}^2)>0$.
Therefore, $\det A =M_1+M_2 > 0$. Since
$\proj_f^-(a_{-i})=-\proj_f^-(a_i)$ for $i\in\Z$ and
 $(a_{-1},a_{-2},a_{-3})=(e,f^{\tau_e},e^{\tau_f\tau_e})$,
 $\proj_f^-(e)$, $\proj_f^-(f^{\tau_e})$ and $\proj_f^-(e^{\tau_f\tau_e}) $
 are linearly independent.
\end{proof}

By Proposition \ref{proj_f^-} and Lemma \ref{independent},
 $|e^T\cup f^T|\leq 6$.
Therefore, by Lemma \ref{orbit}(3), we have the following theorem.
\begin{thm} Let $e$ and $f$ be any Ising vectors of $V$. Then $|\tau_e\tau_f|\leq 6$.
\end{thm}

\subsection{Inner Product}

In this subsection,
we determine the inner product of two Ising vectors.

\begin{thm}
Let $N=|e^T\cup f^T|$. \\
{\rm (1)} If $N=2$, then $\la e,f \ra=0$ or $1/2^5$. \\
{\rm (2)} If $N=3$, then $\la e,f \ra=13/2^{10}$ or $1/2^8$. \\
{\rm (3)} If $N=4$, then $(\la e,f \ra , \la e, e^{\tau_f}\ra
)=(1/2^7,0)$
 or $(1/2^8,1/2^5 )$. \\
{\rm (4)} If $N=5$, then $\la e,f \ra=\la e, e^{\tau_f}\ra=3/2^{9}$. \\
{\rm (5)} If $N=6$, then $\la e,f \ra=5/2^{10}$,
 $\la e,e^{\tau_f}\ra=13/2^{10}$ and
 $\la e^{\tau_f},f^{\tau_e}\ra=1/2^5$.
\end{thm}

\begin{proof}
Let $\lam_1=(e|f)$ and $\lam_2=(e|e^{\tau_f})$ as above. \\
(1) By $N=2$, $|e^T|=|f^T|=1$, that is, $e^{\tau_f}=e$ and $f^{\tau_e}=f$. 
Thus by (\ref{e-f}) and $\lam_2=1$,
\begin{eqnarray*}
0 & = & \frac{1}{7}\left( 2^{11}\lam_1^2 -9\cdot 2^{4}\lam_1 +2^{3}
+\frac{33}{2^4}\right) (e-f)
 + (2^4\lam_1 -\frac{3}{8})(f-e) \\ &  &
 +\frac{1}{2^{4}}(e -f )
 -(2-\frac{1}{8})(e-f) \\
 & = & \frac{1}{7}\left( 2^{11}\lam_1^2 -\cdot 2^{8}\lam_1
\right) (e-f).
\end{eqnarray*}
Hence, $\la e,f \ra=\dfrac{1}{4}\lam_1 = \dfrac{1}{2^5}$ or $0$. \\
(2) Since $f=e^{\tau_f\tau_e}$,
\begin{eqnarray*}
\proj_f^-(f^{\tau_e})=\proj_f^-(e^{\tau_f})=-\proj_f^-(e), \\
\proj_f^-(e^{\tau_f\tau_e})=0 \mbox{ and } \lam_1=\lam_2.
\end{eqnarray*}
Thus by (\ref{proj_f^-}),
\begin{eqnarray*}
0 & = & \frac{1}{7}\left( 2^{15}\lam_1^2 -2^{9}\lam_1 +2^{7}\lam_1
-9 \right) \proj_f^-(e)
 -(2^{8}\lam_1 -5)\proj_f^-(e) \\
 & = & \frac{1}{7}\left( 2^{15}\lam_1^2 -17\cdot 2^{7}\lam_1 +26 \right)
 \proj_f^-(e) \\
  & = & \frac{2}{7}( 2^{8}\lam_1-13)( 2^{6}\lam_1-1).
\end{eqnarray*}
Hence, $\la e,f \ra=\dfrac{1}{4}\lam_1 = \dfrac{13}{2^{10}}$ or $\dfrac{1}{2^{8}}$. \\
(3) Since $|e^T|=|f^T|=2$,
 $e^{\tau_f\tau_e}=e^{\tau_f}$ and $f^{\tau_e\tau_f}=f^{\tau_e}$.
Thus $\proj_f^-(f^{\tau_e})=0$ and
$\proj_f^-(e^{\tau_f\tau_e})=-\proj_f^-(e)$. 
Hence, by \ref{proj_f^-}),
\begin{eqnarray*}
0 & = & \frac{1}{7}\left( 2^{15}\lam_1^2 -2^{9}\lam_1 +2^{7}\lam_2
-9 \right) \proj_f^-(e)
 -\proj_f^-(e) \\
 & = & \frac{2^4}{7}\left( 2^{11}\lam_1^2 -2^{5}\lam_1 +2^{3}\lam_2
-1 \right) \proj_f^-(e).
\end{eqnarray*}
Since $(e^{\tau_f})^{\tau_e}=e^{\tau_f}$ and $e\not= e^{\tau_f}$,
 $\la e, e^{\tau_f}\ra=\dfrac{1}{4}\lam_2=0$ or $\dfrac{1}{2^5}$ by (1).
If $\lam_2=0$,
$$
0 = 2^{11}\lam_1^2 -2^{5}\lam_1 -1 = (2^{5}\lam_1-1)(2^{6}\lam_1+1)
$$
 and so $\la e, f\ra=\dfrac{1}{4}\lam_1 = \dfrac{1}{2^{7}}$ since $\lam_1> 0$.
If $\lam_2=\dfrac{1}{8}$, then $ 0 = 2^{11}\lam_1^2 -2^{5}\lam_1$
 and so $\la e, f\ra=\dfrac{1}{4}\lam_1 = \dfrac{1}{2^{8}}$ since $\lam_1\not= 0$
 by $e\not= e^{\tau_f}$. \\
(4) Since $e^{\tau_f\tau_e}=f^{\tau_e\tau_f}$,
$\proj_f^-(e^{\tau_f\tau_e})=-\proj_f^-(f^{\tau_e})$. Thus by
(\ref{proj_f^-})
\begin{eqnarray*}
0 & = & \frac{1}{7}\left( 2^{15}\lam_1^2 -2^{9}\lam_1 +2^{7}\lam_2
-9 \right) \proj_f^-(e)
 +(2^{8}\lam_1 -5)\proj_f^-(f^{\tau_e})
 -\proj_f^-(f^{\tau_e}) \\
 & = & \frac{1}{7}\left( 2^{15}\lam_1^2 -2^{9}\lam_1 +2^{7}\lam_2
-9 \right) \proj_f^-(e)
 +(2^{8}\lam_1 -6)\proj_f^-(f^{\tau_e})
\end{eqnarray*}
 and so $\la e, f\ra=\dfrac{1}{4}\lam_1 = \dfrac{3}{2^{9}}$
 and $\la e, e^{\tau_f}\ra=\dfrac{1}{4}\lam_2
  =\dfrac{1}{2^{9}}( 9-2^{15}\lam_1^2 +2^{9}\lam_1 )=\dfrac{3}{2^{9}}$. \\
(5) Since $|e^T|=|f^T|=3$, $e^{\rho^{-1}}=e^{\rho^{2}}$ and so
$\proj_f^-(e^{\tau_f\tau_e})=0$. Thus by (\ref{proj_f^-})
\begin{eqnarray*}
0 & = & \frac{1}{7}\left( 2^{15}\lam_1^2 -2^{9}\lam_1 +2^{7}\lam_2
-9 \right) \proj_f^-(e)
 +(2^{8}\lam_1 -5)\proj_f^-(f^{\tau_e})
\end{eqnarray*}
and so $\la e, f\ra=\dfrac{1}{4}\lam_1 = \dfrac{5}{2^{10}}$ and $\la
e, e^{\tau_f}\ra=\dfrac{1}{4}\lam_2=\dfrac{1}{2^{9}}(
9-2^{15}\lam_1^2 +2^{9}\lam_1 )=\dfrac{13}{2^{10}}$. By
(\ref{(a|alpha(a,x))}) and (\ref{(e|alpha(f,f^{tau_e}))}),
\begin{eqnarray*}
(e|\alpha(e,e^{\tau_f})) & = & \frac{31}{16}\lam_2-\frac{1}{16} \\
(e|\alpha(f,f^{\tau_e})) & = & \frac{15}{8}\lam_1
-\frac{1}{16}(1-\lam_2)
\end{eqnarray*}
and so by (\ref{e-f}),
\begin{eqnarray*}
( e | f^{\tau_e\tau_f} ) & = & \lam_2
 +\frac{2^4}{7}(2^{11}\lam_1^2 -9 \cdot 2^4\lam_1 +\frac{33}{2^{4}}
 +2^3\lam_2 ) (1-\lam_1) \\
 &  & +2^8(\lam_1-\frac{3}{2^7})(\lam_1-\lam_2)
 +2^4(\frac{15}{8}\lam_1 -\frac{1}{16}(1-\lam_2)
  -\frac{31}{16}\lam_2+\frac{1}{16}) \\
 & = & \frac{1}{8}.
\end{eqnarray*}
Thus, $\la e^{\tau_f},f^{\tau_e}\ra
 =\la e,f^{\tau_e\tau_f}\ra=\dfrac{1}{2^5}$.
\end{proof}

\begin{rmk}
(1) and (2) of the last Theorem was shown by Miyamoto
 in \cite{M1} and \cite{M2}.
\end{rmk}

\end{document}